\documentclass[12pt]{article}
\usepackage{amssymb}
\usepackage{amsmath}
\usepackage{latexsym}
\usepackage{mathrsfs}
\usepackage{color}

\numberwithin{equation}{section}
\newtheorem{thm}{Theorem}

\newtheorem{lem}[thm]{Lemma}

\newcommand{\A}{{\mathcal{A}}}

\newcommand{\R}{{\mathbb{R}}}
\newcommand{\N}{{\mathbb{N}}}

\def\Int{\displaystyle\int}

\newcommand{\beq}{\begin{equation}}
\newcommand{\eeq}{\end{equation}}

\title{Monotonicity, continuity and differentiability results for the  $L^p$ Hardy constant}

\author{
Gerassimos Barbatis\footnote{Department of Mathematics,
 University of Athens,  15784 Athens, Greece}\ \ 
and Pier Domenico Lamberti 
 \footnote{Dipartimento di Matematica,
Universit\`{a} degli Studi di Padova,
Via Trieste 63, 35121 Padova, Italy}}

\date{\ }

\begin{document}

\newcommand{\rea}{\mathbb{R}}

\maketitle

%\subjclass{Primary 30C20, 47H30, 45F15; Secondary 35J25}
%
%\date{September 12, 2001}
%
%\thanks{}
%

\noindent
%{\bf Keywords:}

%\vspace{6pt}
%\noindent
%{\bf 2000 Mathematics Subject Classification:}
% 35P15, 35J40, 47A75, 47B25.
%
%

\begin{abstract}
\noindent We consider the $L^p$ Hardy inequality involving the distance to the boundary for a domain in the $n$-dimensional Euclidean space. We study the dependence
on $p$ of the corresponding best constant and we prove monotonicity, continuity and differentiability results. The focus is on non-convex domains in which case such constant is in general not  explicitly known.   
\end{abstract}

\vspace{11pt}


\vspace{6pt}

\parskip=5pt

\section{Introduction}

Given a bounded domain $\Omega $ in ${\mathbb{R}}^n$ and 
$p\in ]1,\infty [$, we say that the $L^p$  Hardy inequality holds in $\Omega$
if there exists $c>0$ such that 
\begin{equation}
\label{mainhardy}
\int_{\Omega }|\nabla u|^pdx \geq c \int_{\Omega }\frac{|u|^p}{d^p}dx \, , \quad \mbox{ for all $u\in C^{\infty }_c(\Omega)$}, 
\end{equation}
where $d(x)= {\rm dist }(x, \partial \Omega)$, $x\in \Omega$. The $L^p$ Hardy constant of $\Omega$ is 
the best constant for inequality (\ref{mainhardy}) and is denoted here by  $H_p$.

It is well-known that the $L^p$ Hardy inequality holds for all $p \in ]1,\infty[$ under  weak regularity assumptions on $\Omega$, for example if $\Omega$ has a Lipschitz boundary. Moreover, if $\Omega$ is convex, and more generally if it is weakly mean convex, i.e. if $\Delta d\le 0$ in the distributional sense in $\Omega$, then $H_p=((p-1)/p)^p$; see \cite{MMP,barfilter}. 
If $\Omega$ is not weakly mean convex, little is known about the precise value of $H_p$ and the available results   only hold  for $p=2$ and for 
special domains, for example  circular sectors and quadrilaterals in the plane. We refer to  \cite{ avk, barlam, barfilter, BarTer,   BarTer2, bremar, davies, laso, MMP} for more information. 
We also refer to the monograph \cite{permalkuf} for an introduction to the study of Hardy and Hardy-type inequalities with a historical perspective.

In this article we study the dependence of $H_p$ upon variation of $p$ and we prove four main results. First, we prove  that 
$p(1+H_p^{1/p})$ is a non decreasing function of  $p\in ]1,\infty [$, and this is done without any smoothness assumption on $\Omega$, see Theorem~\ref{thm:mono}. In particular, it easily follows that $H_p$ is right-continuous at any point  $p\in ]1,\infty[$.  Second,
we prove that   if $\Omega $ is of class $C^2$ then $H_p$ is also left-continuous, hence it is continuous  on $ ]1,\infty [$, see Theorem~\ref{cont}. Third, we prove that if $\Omega $ is of class $C^2$ then $H_p$ is differentiable at any point $p\in ]1,\infty [$ such that $H_p< ((p-1)/p)^p$, and we compute a formula for the corresponding derivative, see Theorem~\ref{diff}. 

We note that the proofs of our continuity and differentiability results exploit a result by \cite{MMP}, where it was shown in particular that if $H_p<((p-1)/p)^p$ then equality  is attained in (\ref{mainhardy})  for some function $u_p\in W^{1,p}_0(\Omega)$  which behaves like $d^{\alpha}_{\Omega}$ near $\partial\Omega$ for a suitable $\alpha \in ]0,1[$.  Importantly, the results of \cite{MMP} are proved under the assumption that  $\Omega$ is of class $C^2$, and removing that assumption is not easy.  
The function $u_p$ is uniquely identified by the extra normalizing conditions $u_p>0$ and $\int_{\Omega}u_p^p/d^pdx=1$. The fourth main result of the paper is a continuity result
for the dependence of $u_p$ and $\nabla u_p$ on $p$, see Theorem~\ref{conte}.

As is well-known, if  equality  is attained in (\ref{mainhardy}) for some nontrivial function $u\in W^{1,p}_0(\Omega)$, then $u$ is a minimizer for the  Hardy quotient
\begin{equation}
\label{hardyquo}
R_p[u]:=\frac{\int_{\Omega}|\nabla u|^pdx}{\int_{\Omega}\frac{|u|^p}{d^p}dx } 
\end{equation}
and solves the equation

\begin{equation}\label{pde}
-\Delta_p u=H_p\frac{|u|^{p-2}u}{d^p},
\end{equation}
where $\Delta _pu={\rm div}(|\nabla u|^{p-2}\nabla u)$ is the $p$-Laplacian.

Problem (\ref{pde}) is a singular variant of the well-known  eigenvalue problem for the Dirichlet $p$-Laplacian 
\begin{equation}\label{pdeplap}
-\Delta_p u=\lambda_p|u|^{p-2}u,
\end{equation}
where $H_p$ is replaced by the first eigenvalue $\lambda_p$ of the $p$-Laplacian, which in turn is the minimum over
$W^{1,p}_0(\Omega)\setminus\{0\}$ of the Rayleigh quotient
\begin{equation}
\label{rayplap}
\frac{\int_{\Omega}|\nabla u|^pdx}{\int_{\Omega}|u|^pdx } .\end{equation}

The study of the dependence of $\lambda_p$ on $p$ was initiated in the article \cite{lqv} which has inspired many authors, ourselves included.    We refer to \cite{anello, degio, erc} for recent closely related results.
In fact, the proofs of our monotonicity and continuity results exploit some ideas of \cite{lqv}. However, we point out that although the two problems (\ref{pde}) and (\ref{pdeplap}) look similar, they  are radically different. 
For example, if  $\Omega$ has finite Lebesgue measure, the Rayleigh quotient (\ref{rayplap}) has always a minimizer and if  $\Omega$ is also sufficiently smooth,  the gradient of such minimizer does not blow up at the boundary. As is well-known, one of the main differences between the two problems is related to the lack of compactness for the embedding of the Sobolev space $W^{1,p}_0(\Omega)$ into the natural weighted  space $L^p(\Omega, d^{-p}dx)$, which is also responsible for the appearence of a large essential spectrum for problem (\ref{pde}) in the case $p=2$.  Thus, the study of the dependence of $H_p$  on $p$, leads to a number of difficulties which require a detailed analysis. 

We point out the our differentiability result can also be proved, with obvious simplifications, for the dependence of $\lambda_p$ on $p$. Since we have not 
found such result in the literature, we find it natural to state it in the Appendix.

\section{Preliminaries}

Unless otherwise indicated, by $\Omega$ we  denote a bounded domain (i.e. a bounded open connected set)  in $\R^n$.
If $p\in ]1,+\infty[$ we denote by $W^{1,p}(\Omega)$ the standard Sobolev space and by $W^{1,p}_0(\Omega)$ the closure in  $W^{1,p}(\Omega)$ of the space $C^{\infty}_c(\Omega )$ of all $C^{\infty}$-functions with compact support in $\Omega$.

The $L^p$ Hardy constant is defined by
\begin{equation}
H_p =\inf_{u\in W^{1,p}_0(\Omega), u\neq 0} R_{p}[u],
\label{ray}
\end{equation}
and  if $H_p>0$ we say that the $L^p$ Hardy inequality is valid on
$\Omega$.  

It is well known that if $\Omega$ has a Lipschitz continuous boundary then $0<H_p\leq ((p-1)/p)^p$. It is also known that if  $\Omega$ is of class $C^2$ then there exists a minimizer $u$ in (\ref{ray}) if and only if $H_p< ((p-1)/p)^p$,  see \cite{MMP,mash}; moreover,  the minimizer is unique up to a multiplicative constant, can be chosen to be positive  
 and there exists $c>0$ such that
\begin{equation}\label{decay}
c^{-1}d(x)^{\alpha_p} \leq u(x) \leq  c d(x)^{\alpha_p},\ \ x\in \Omega ,
\end{equation}
where $\alpha_p\in \, ](p-1)/p , 1[$  denotes the largest solution to the equation 
\begin{equation}
(p-1)\alpha^{p-1}(1-\alpha)=H_p.
\label{alpha}
\end{equation}
We set for simplicity
\[
\A=\left\{p\in ]1,\infty [:\ H_p<((p-1)/p)^p \right\}\, .
\]
In the sequel and provided $\Omega$ is $C^2$ we shall denote for any $p\in\A$ by $u_p$ the positive minimizer normalized by the condition $\int_{\Omega}|u_p/d|^pdx=1$.  Inequalities (\ref{decay}) suggest that 
$\nabla u_p$ behaves like $d^{\alpha_p-1}$ close to the boundary of $\Omega$. In fact we can prove the following
lemma which is a variant of \cite[Thm.~4]{barlam} providing further information on the dependence of the constants on $p$. We emphasize that in this lemma we do not assume that $H_p$ depends continuously on $p$.
\begin{lem}\label{lemgrowth} Assume that $\Omega$ is of class $C^2$ and $p_0\in\A$. There exists $c>0$ such that 
\begin{equation}
\label{lemgrowth0} 
u_p(x)\le c d^{\alpha_p}(x),\ \ \  |\nabla u_p(x)| \le c d^{\alpha_p-1}(x),
\end{equation}
for all $p\in \A$ sufficiently close to $p_0$ and for all $x\in \Omega$.
In particular, $u_p\in W^{1,q}_0(\Omega)$ for all $q\in [1, 1/(1-\alpha_p)[$.
\end{lem}
{\em Proof.} The existence for each $p\in\A$ of a constant $c=c(p)>0$ such that the first inequality in (\ref{lemgrowth0})
holds has been proved in
\cite[Lemma 9]{MMP} and \cite[Lemma 5.2]{mash}. The existence for each $p\in\A$ of a constant $c=c(p)>0$ such that the second inequality in (\ref{lemgrowth0}) holds has been proved in \cite[Theorem 4]{barlam}. We shall now show that $c(p)$ can be chosen so that it is locally bounded with respect to $p\in\A$.

Let $p\in\A$ and let $u\in W^{1,p}_0(\Omega)$ be a positive minimizer of the $L^p$ Hardy constant normalized by $\int_{\Omega}u^p/d^pdx=1$. Let $\alpha$ be as in (\ref{alpha}). For any $\beta>0 $, we set $\Omega_{\beta}=\{x\in \Omega :\ d(x)<\beta\}$. Let $\beta_0>0$ be small enough so that $d(x)$ is twice continuously differentiable in $\Omega_{2\beta_0}$. Following \cite{MMP, mash},  we define
\[
v=d^{\alpha} (1-d). 
\]
A direct computation gives that in $\Omega_{2\beta_0}$,
\begin{eqnarray}
&&\hspace{-1.5cm} -\Delta_p v - H_p \frac{v^{p-1}}{d^p} =\nonumber \\[0.2cm]
&=& (p-1)\alpha^{p-1}d^{\alpha p-\alpha-p}  \left\{ (1-\alpha) \left[ \bigg(    1-(1+\frac{1}{\alpha})d \bigg)^{p-1} -\Big(1-d\Big)^{p-1}  \right] 
\right. \nonumber \\[0.2cm]
&& \left.  + (1+\frac{1}{\alpha}) \bigg(    1-(1+\frac{1}{\alpha})d\bigg)^{p-2}d   \right\}  \nonumber \\[0.2cm]
&& -\alpha^{p-1}d^{\alpha  p -\alpha -p +1}\bigg(    1-(1+\frac{1}{\alpha})d\bigg)^{p-1} \Delta d  \nonumber \\[0.2cm]
&=&  d^{\alpha p-\alpha-p}   (A+ B \, d\Delta d)  \, ,  \label{2hm}
\end{eqnarray}
where terms in $A$ do not involve $\Delta d$. We expand $A$ in powers of $d$ and obtain
\begin{eqnarray*}
A&=& (p-1)\alpha^{p-2}  ( \alpha p -p+2 ) d +O(d^{2 })  \\[0.2cm]
&\geq& (p-1)\alpha^{p-2}  d + O(d^{2 }).
\end{eqnarray*}
It can easily be verified that the coefficient of $d^{2}$ is locally bounded with respect to $p\in ]1,+\infty[$. Hence there exists $\beta_1\in ]0,\beta_0[$ which is locally bounded away from zero with respect to $p$ such that
\begin{equation}
A \geq  \frac{(p-1)\alpha^{p-2} }{2}d \; , \qquad \mbox{ in }\Omega_{\beta_1}.
\label{bb}
\end{equation}
Since $\Delta d$ is bounded in $\Omega_{\beta_0}$, it follows from (\ref{2hm}) and (\ref{bb}) that there exists
$\beta_2\in ]0,\beta_1[$  bounded away from zero locally in $p\in \A$
such that 
\[
-\Delta_p v - H_p \frac{v^{p-1}}{d^p} \geq 0 \;\; , \;\; \mbox{ in }\Omega_{\beta_2}.
\]
Now, let
\[
C_1(p) = \sup\big\{ u(x) \, : \;  x\in \{d(x)=\beta_2\} \big\}\, .
\]
The constant $C_1(p)$ is finite by standard regularity results for quasilinear elliptic equations.
Looking e.g. at the proof of Theorems 1 and 2  of the classical paper of Serrin \cite{se} we can trace the dependence of $C_1(p)$ in $p$ for $p\le n$ and see that it is locally bounded for  $p\le n$.
As mentioned in \cite{se}, the case $p>n$ is simpler since the result follows by the Sobolev embedding. We note that  the fact that the Sobolev constant blows-up as $p \to n^+$ is not a problem, since the argument used in \cite[Theorem 2]{se} for $p=n$ can be extended without changes to include all $p$ in a neighborhood of $n$. We omit the details.

Defining next $C^*=C_1/(\beta_2^{\alpha} (1-\beta_2))$, we then have
\[
C^*=\sup \Big\{ \frac{u(x)}{v(x)} \, , \;  x\in \{d(x)=\beta_2\} \Big\}. 
\]
Applying \cite[Proposition 3.1]{mash} we conclude that
\[
u(x)\leq C^*v(x) \leq C^* d^{\alpha} \; , \quad \mbox{ in }\Omega_{\beta_2}.
\]
This estimate clearly holds true also in $\Omega\setminus\Omega_{\beta_2}$, with a constant $C^*$ still remaining locally bounded with respect to $p\in \A$, completing the proof of the first estimate of (\ref{lemgrowth0}).

For the second inequality we apply the regularity estimates of \cite[Theorems 1.1 and 1.2]{dumi1}, as was done in
\cite{barlam}. The constants involved are locally bounded in $p$ (see in particular \cite[Remark 5.1]{dumi1}). This completes the proof.
$\hfill\Box$

\section{Monotonicity and continuity of the Hardy constant}

The following theorem holds without any smoothness assumption of $\Omega$ (not even the boundedness of $\Omega$ is actually required) and is inspired by the monotonicity result proved in Lindqvist~\cite{lqv} for the 
first eigenvalue of the $p$-Laplacian. 

\begin{thm}
\label{thm:mono}
The function
\[
p\mapsto p(1+H_p^{1/p})
\]
is non-decreasing in $]1,+\infty[$.
\end{thm}
{\em Proof.} Let $1<p<s$ and let $\psi\in C^{\infty}_c(\Omega)$. Then the function
\[
u =|\psi|^{\frac{s}{p}}d^{1-\frac{s}{p}}
\]
belongs to $W^{1,p}_0(\Omega)$ and
\begin{eqnarray*}
\Big(\int_{\Omega}|\nabla u|^pdx\Big)^{1/p} &=& \Big(\int_{\Omega}\Big| \frac{s}{p}\big(\frac{|\psi|}{d}\big)^{\frac{s}{p}-1}\nabla\psi  + (1-\frac{s}{p})\big(\frac{|\psi|}{d}\big)^{\frac{s}{p}}\nabla d \Big|^pdx \Big)^{1/p}\\
&\leq& \frac{s}{p} \Big(\int_{\Omega}\big(\frac{|\psi|}{d}\big)^{s-p}|\nabla\psi|^p dx\Big)^{1/p} + 
\frac{s-p}{p}\Big(\int_{\Omega}\big(\frac{|\psi|}{d}\big)^{s} dx \Big)^{1/p}\\
&\leq& \frac{s}{p} \Big(\int_{\Omega}|\nabla\psi|^s dx\Big)^{1/s}
\Big(\int_{\Omega}\big(\frac{|\psi|}{d}\big)^{s} dx \Big)^{\frac{1}{p}-\frac{1}{s}}
+  \frac{s-p}{p}\Big(\int_{\Omega}\big(\frac{|\psi|}{d}\big)^{s} dx \Big)^{1/p}.
\end{eqnarray*}
This implies
\[
H_p^{1/p}\leq R_p[u]^{1/p}\leq \frac{s}{p}R_s[\psi]^{1/s} +\frac{s-p}{p}.
\]
Taking the infimum over all $\psi\in C^{\infty}_c(\Omega)$ we conclude that
\[
H_p^{1/p}\leq \frac{s}{p}H_s^{1/s} +\frac{s-p}{p},
\]
and the result follows. \hfill$\Box$ 

\noindent{\bf Remarks.} (1) For $\alpha \in [0,1]$ let
\[
\lambda_{\alpha ,p} = \inf_{u\in W^{1,p}_0(\Omega), u\neq 0}\frac{\Int_{\Omega}|\nabla u|^pdx}{\Int_{\Omega}
\frac{|u|^p}{d^{ap}}dx};
\]
so $\lambda_{1,p}=H_p$ and $\lambda_{0,p}=\lambda_p$ is the first eigenvalue of the Dirichlet $p$-Laplacian in $\Omega$ (see the introduction).
It has been shown in \cite[Theorem 3.2]{lqv} that the function $p\mapsto p\lambda_{0,p}^{1/p}$ is non-decreasing in $]1,\infty[$. In view of this and Theorem \ref{thm:mono} it is tempting to believe that for any fixed $\alpha\in [0,1]$ the map
$p\mapsto  p(\alpha+\lambda_{\alpha,p}^{1/p})$ is  non-decreasing in $]1,\infty[$. However it can be seen that the method of proof fails
for $\alpha\in ]0,1[$.

\noindent (2) It follows from Theorem \ref{thm:mono} that the function $p\mapsto H_p$ has one-sided limits at every $p>1$ and
\begin{equation}
\lim_{s\to p-}H_s \leq H_p \leq \lim_{s\to p+}H_s \; .
\label{ra}
\end{equation}

\begin{lem}
\label{lem1}
We have
\[
\limsup_{s\to p}H_s = \lim_{s\to p+}H_s = H_p \; .
\]
\end{lem}
{\em Proof.} Given any $u\in C^{\infty}_c(\Omega)$ we have $H_s\leq R_s[u]$ and therefore
\[
\limsup_{s\to p}H_s \leq R_p[u]. 
\]
Taking the infimum over all $u\in C^{\infty}_c(\Omega)$ we obtain $\limsup_{s\to p}H_s \leq H_p$ which combined with (\ref{ra})
yields the result. $\hfill\Box$

In order to prove Theorem~\ref{cont} we need the following lemmas. The first can be proved simply by differentiating under the integral sign.

\begin{lem}
\label{prop1}
Let $u\in W^{1,p}_0(\Omega)$ be fixed. The functions defined by
\[
N(s) =\int_{\Omega}|\nabla u|^s dx \quad , \qquad D(s)=\int_{\Omega}\frac{|u|^s}{d^s}dx
\]
are differentiable in $]1,p[$ and
\[
 N'(s)=\int_{\Omega} |\nabla u|^s \ln |\nabla u| dx \quad , \qquad
 D'(s)= \int_{\Omega} \frac{|u|^s}{d^s} \ln \big(\frac{|u|}{d}\big) dx
\]
for all $1<s<p$.
\end{lem}

\begin{lem}
\label{lem2}
Assume that $\Omega$ is of class $C^2$. We have
\[
\liminf_{s\to p}H_s \geq H_p\, .
\]
\end{lem}
{\em Proof.} It follows from (\ref{ra}) that
\[
\liminf_{s\to p}H_s =\liminf_{s\to p-}H_s \; .
\]
Suppose by contradiction that this liminf is a number $L<H_p\; $.
Let $s_n$, $n\in \N$,  be an increasing sequence of exponents with $s_n\to p$ and $H_{s_n}\to L$ as $n\to\infty$. Then, since $L<H_p\le (\frac{p-1}{p})^{p}$, we have that
$H_{s_n}<(\frac{s_n-1}{s_n})^{s_n}$ for all $n\in\N$ sufficently large  and therefore the $L^{s_n}$-Hardy quotient has a positive minimizer $u_{s_n}$. Let $\alpha_{s_n}$ be the corresponding exponents defined as in (\ref{alpha}). It then follows that
$\lim_{n\to \infty }\alpha_{s_n} >(p-1)/p \, $. 
Applying Lemma~\ref{lemgrowth} we thus obtain that
\begin{equation}
\|u_{s_n}\|_{W^{1,p+\epsilon}_0(\Omega)}\leq M
\label{1}
\end{equation}
for some fixed $\epsilon,M>0$ and all $n\in\N$ sufficiently large. Hence
\begin{eqnarray*}
H_p &\leq& \liminf_{n\to \infty } R_p[u_{s_n}] \\
&=& \liminf_{n\to \infty } \Big(  R_{s_n}[u_{s_n}]  +  \big\{ R_p[u_{s_n}]-R_{s_n}[u_{s_n}]  \big\} \Big) \\
&=& L+  \liminf_{n\to \infty } \big(  R_p[u_{s_n}]-R_{s_n}[u_{s_n}] \big).
\end{eqnarray*}
To reach a contradiction it is enough to prove that the last liminf is zero. Now, by Lemma \ref{prop1} and (\ref{1}) the function $s\mapsto R_s[u_{s_n}]$ is differentiable in $(s_n,p)$ for each fixed $n\in\N$. Hence by the Mean Value Theorem,
for each $n\in\N$ there exists $\xi_n\in (s_n,p)$ such that
\[
R_p[u_{s_n}]-R_{s_n}[u_{s_n}] = (p-s_n)\frac{dR_p[u_{s_n}]}{dp}\bigg|_{p=\xi_n} \; .
\]
From Lemma \ref{prop1} and (\ref{1}) easily follows that $\frac{dR_p[u_{s_n}]}{dp}\bigg|_{p=\xi_n}$ remains bounded as $n\to\infty$. This concludes the proof. $\hfill\Box$

\begin{thm}
\label{cont}
Let $\Omega$ be bounded with $C^2$ boundary. Then the function $p\mapsto H_p$ is continuous on $]1,\infty[$.
\end{thm}
{\em Proof.} Follows from Lemmas \ref{lem1} and \ref{lem2}. $\hfill\Box$

\section{Differentiability of the Hardy constant}

We recall that $\A=\left\{p\in ]1,\infty [ \; : \; H_p<((p-1)/p)^p \right\}$.
The proof of the following theorem is based on adapting the arguments of Lindqvist~\cite[Thm.~3.6]{lqv}. 
\begin{thm}\label{conte}Let $\Omega$ be of class $C^2$ and  $p_0\in \A$. Then for all $p$ sufficiently close to $p_0$ we have $p\in\A$ and
$u_p,u_{p_0}\in W^{1,\max\{p_0,p\}}(\Omega )$.  Moreover
\begin{equation}\label{conte0}
\lim_{p\to p_0}\| u_p-u_{p_0}\|_{W^{1,\max\{p_0,p\}}(\Omega )}=0.
\end{equation} 
\end{thm}
{\em Proof.} Theorem~\ref{cont} and Lemma~\ref{lemgrowth} easily imply that
for $p$ close enough to $p_0$ we have $p\in\A$ and, moreover, $u_p\in W^{1,p_0}(\Omega)$ and $u_{p_0}\in W^{1,p}(\Omega)$.

We now prove (\ref{conte0}). 
Let $\delta >0$ be fixed in such a way that $p_0+2\delta <1/(1-\alpha_{p_0})$. By Theorem~\ref{cont} and Lemma~\ref{lemgrowth} it follows that there exists a constant $c>0$ independent of $p$ such 
that
\[
\| u_p\|_{W^{1, p_0+\delta }(\Omega)}\le c \, ,
\]
for all $p\in \A$ sufficiently close to $p_0$. Moreover, since $\Omega$ has $C^2$ boundary we have
$u_p\in W_0^{1, p_0+\delta }(\Omega)$ for any such $p$. 

By the reflexivity of the space $W^{1, p_0+\delta }_0(\Omega)$ and the Rellich-Kondrachov Theorem it follows that there exists $\tilde u\in W^{1, p_0+\delta }_0(\Omega)$ such that, up to taking a subsequence, $\nabla u_p\rightharpoonup \nabla\tilde u$ weakly in $L^{p_0+\delta}(\Omega)$ and $u_p\to \tilde u$ in $L^{p_0+\delta}(\Omega)$ as $p\to p_0$. Note that $\int_{\Omega}|\tilde u|^{p_0}/d^{p_0}dx=1$, which can be deduced by passing to the limit as $p\to p_0$ in the equality $\int_{\Omega}|u_p|^p/d^pdx=1$ and using the Dominated Convergence Theorem combined with estimates (\ref{lemgrowth0}).  In particular $\tilde u\ne 0$.  Clearly,  $\nabla u_p\rightharpoonup \nabla\tilde u$ weakly in $L^{p_0}(\Omega)$  hence
\begin{equation}\label{conte1}\int_{\Omega}|\nabla \tilde u|^{p_0}dx\le \liminf_{p\to p_0}\int_{\Omega}|\nabla u_p|^{p_0}dx \end{equation}
as $p\to p_0$. By the Mean Value Theorem and Lemma~\ref{prop1} we have that 
\begin{eqnarray}\label{conte2}
\int_{\Omega}|\nabla u_p|^{p_0}dx& =&  \int_{\Omega}|\nabla u_p|^{p}dx+(p_0-p)\int_{\Omega }|\nabla u_p|^{s_p}\ln |\nabla u_p|dx\nonumber  \\
& =& H_p+ (p_0-p)\int_{\Omega }|\nabla u_p|^{s_p}\ln |\nabla u_p|dx ,
\end{eqnarray}
for some real number $s_p$ between $p_0$ and $p$. It is clear that by the uniform  boundedness of the norms of $u_p$ in $W_0^{1,p_0+\delta}(\Omega)$, the  integrals $\int_{\Omega }|\nabla u_p|^{s_p}\ln |\nabla u_p|dx$ are uniformly bounded for $p$ close enough to $p_0$. Thus,  by passing to the limit as $p\to p_0$ in (\ref{conte2}) and using the continuity of the map $p\mapsto H_p$ it follows that
\begin{equation}\label{preclar}
\lim_{p\to p_0}\int_{\Omega}|\nabla u_p|^{p_0}dx= \lim_{p\to p_0}H_p  =H_{p_0}.
\end{equation}
This combined with (\ref{conte1}) and condition   $\int_{\Omega}|\tilde u|^{p_0}/d^{p_0}dx=1$     implies that $\int_{\Omega}|\nabla\tilde u|^{p_0}=H_{p_0}$. Thus, $\tilde u =u_{p_0}$. 

As in \cite[Thm.~3.6]{lqv} we now use Clarkson's inequalities. If $\max\{ p_0, p\}\geq 2$ we have 
\begin{eqnarray}
&& \int_{\Omega}\left| \frac{\nabla u_p-\nabla u_{p_0}}{2}  \right|^{\max\{ p_0, p\}}dx   \nonumber\\ 
&&\qquad \le \frac{1}{2}\int_{\Omega}|\nabla u_p|^{\max\{ p_0, p\}}dx   + \frac{1}{2}\int_{\Omega}|\nabla u_{p_0}|^{\max\{ p_0, p\}}dx
  \nonumber\\ 
&& \qquad
-\int_{\Omega}\left| \frac{\nabla u_p+\nabla u_{p_0}}{2}  \right|^{\max\{ p_0, p\}}dx \le \frac{1}{2}\int_{\Omega}|\nabla u_p|^{\max\{ p_0, p\}}dx  
\nonumber \\
&&\qquad +\frac{1}{2}\int_{\Omega} |\nabla u_{p_0}|^{\max\{ p_0, p\}}dx     -H_{\max\{p,p_0\}}\int_{\Omega}\left|\frac{u_p+u_{p_0}}{2d}\right|^{\max\{ p_0, p\}}dx 
\label{clark}
\end{eqnarray}
By the continuity of the $L^p$-norm , it follows that
\begin{equation}\label{postclar}\lim_{p\to p_0}
\int_{\Omega} |\nabla u_{p_0}|^{\max\{ p_0, p\}}dx= \int_{\Omega} |\nabla u_{p_0}|^{p_0}dx=H_{p_0}.
\end{equation}
Moreover, 
using the Dominated Convergence Theorem combined with estimates (\ref{lemgrowth0}) yields
$$\lim_{p\to p_0}
\int_{\Omega}\left|\frac{u_p+u_{p_0}}{2d}\right|^{\max\{ p_0, p\}}dx = \int_{\Omega}\left|\frac{u_{p_0}}{d}\right|^{p_0}dx=1.
$$
We then deduce from  (\ref{preclar})-(\ref{postclar}) and Theorem~\ref{cont} that $\int_{\Omega}\left| \frac{\nabla u_p-\nabla u_{p_0}}{2}  \right|^{\max\{ p_0, p\}}dx \to 0$ as required. The case $p_0<2$ can be treated in a similar way using the appropriate Clarkson inequality for $p<2$.\hfill $\Box$

\begin{thm}\label{diff}
Let $\Omega$ be of class $C^2$ . Then the map $p\mapsto H_p$ is of class $C^1$ on $\A$ and 
\begin{equation}\label{diff0}
H_p'=\int_{\Omega}|\nabla u_p|^p\ln |\nabla u_p|dx-H_p\int_{\Omega}\frac{u_p^p}{d^p}\ln \frac{u_p}{d^p}dx \, , \quad
p\in\A\, .
\end{equation}
\end{thm}
{\em Proof.}  Let $p_0\in \A$ be fixed. Since $\A$ is an open set, if $p>1$ is sufficiently close to $p_0$, we have that $p\in \A$ hence the minimizer $u_p$ exists.
Moreover, by Lemma~\ref{lemgrowth} and Theorem~\ref{cont}, there exist $\epsilon, \delta>0$ such that   $p< 1/(1-\alpha_{p_0}) +\epsilon$ and
\begin{equation}
\label{diff1}
u_p\in W^{1,1/(1-\alpha_{p_0}) +\epsilon}(\Omega) ,
\end{equation}
for all $p\in ]p_0-\delta, p_0+\delta [$.  
Since $u_{p_0}$ and $u_{p}$ minimize the corresponding Rayleigh quotients, we have
\begin{equation}
\label{diff2}
R_p[u_p]-R_{p_0}[u_{p}]\le H_p-H_{p_0}\le R_p[u_{p_0}]-R_{p_0}[u_{p_0}].
\end{equation} 
By (\ref{diff1}) and Lemma~\ref{prop1} 
we have that for any fixed $p\in ]p_0-\delta , p_0+\delta[$, the maps $q\mapsto R_q[u_p]$
are differentiable on  $]p_0-\delta, p_0+\delta [$, hence (\ref{diff2}) implies that 
\begin{equation}
\label{diff3}
R'_{p_{\xi}}[u_p](p-p_0)\le H_p-H_{p_0}\le R'_{p_{\eta}}[u_{p_0}](p-p_0)
\end{equation} 
for some $p_{\xi},p_{\eta}$ between $p_0$ and $p$. By Theorem~\ref{conte} and estimates (\ref{lemgrowth0}) one can prove that 
\begin{equation}\label{diff3,5}
R'_{p_{\xi}}[u_p],\  R'_{p_{\eta}}[u_{p_0}]\to  R'_{p_{0}}[u_{p_0}],\ \ {\rm as}\ p\to p_0.
\end{equation}
Indeed, by (\ref{conte0}) it follows that possibly passing to  subsequences $\lim_{p\to p_0}u_p(x)=u_{p_0}(x)$ a.e. in $\Omega$ 
which combined with estimates (\ref{lemgrowth0}) allows passing to the limit under the integral signs in order to get (\ref{diff3,5}).  Thus,  (\ref{diff3}) and (\ref{diff3,5})  imply that $H_p$ is differentiable at $p=p_0$.  Formula (\ref{diff0}) for  $p=p_0$ is then easily proved by using the 
formulas provided by Lemma~\ref{prop1}. 

Finally, in order to prove that the map $p\mapsto H'_p$ is continuous on $\A$, one has simply to apply again   Theorem~\ref{conte} combined with estimates (\ref{lemgrowth0}) as above. \hfill $\Box$

{\noindent\bf Remarks.} 
(1) We note explicitly that since $H_p=\int_{\Omega}|\nabla u_p|^pdx$ we have that 
\begin{eqnarray}
\int_{\Omega}|\nabla ku_p|^p\ln |\nabla ku_p|dx-H_p\int_{\Omega}\frac{|ku_p|^p}{d^p}\ln \frac{|ku_p|}{d}dx\nonumber \\
=|k|^p\left(\int_{\Omega}|\nabla u_p|^p\ln |\nabla u_p|dx-H_p\int_{\Omega}\frac{|u_p|^p}{d^p}\ln \frac{|u_p|}{d}dx\right)
\end{eqnarray}
for any $k\in {\mathbb{R}}$, with $k\ne 0$. In particular, it follows that if we consider a minimizer $u$ for $H_p$ which is not necessarily normalized as $u_p$ then 
\begin{equation}\label{remder}
H'_p=\frac{\int_{\Omega}|\nabla u|^p\ln |\nabla u|dx}{\int_{\Omega}\frac{|u|^p}{d^p}dx}-\frac{H_p\int_{\Omega}\frac{|u|^p}{d^p}\ln \frac{|u|}{d^p}dx}{\int_{\Omega}\frac{|u|^p}{d^p}dx}.
\end{equation}

(2)  For all $p\in \A $ any minimizer $u$ for $H_p$ satisfies the following inequality

\begin{eqnarray}
\label{maincor0}
\lefteqn{H_p
\int_\Omega \frac{|u|^p}{d^p}
\ln \frac{|u|}{d}\,dx
}\nonumber\\ & & 
\le
\left(
H_p+H_p^{\frac{p-1}{p}}
-\frac{H_p\ln H_p}{p}
\right)
\int_\Omega \frac{|u|^p}{d^p}\,dx +\int_\Omega |\nabla u|^p\ln|\nabla u|\,dx.
\end{eqnarray}

Indeed, by  Theorems~\ref{thm:mono}, \ref{diff} the derivative of the function $p\mapsto p(1+H_p^{1/p})$ is non-negative, hence inequality (\ref{maincor0}) follows by
formula (\ref{remder}).

\section{Appendix}

The proof of Theorem~\ref{diff} can be carried out also in the case of the first eigenvalue $\lambda_p$ of the $p$-Laplacian defined by
\begin{equation}
\label{rayplapapp}\lambda_p=\inf_{v\in W^{1,p}_0(\Omega), \ v\ne 0 }
\frac{\int_{\Omega}|\nabla v|^pdx}{\int_{\Omega}|v|^pdx },\end{equation}
see the Introduction. Recall that if $\Omega $ is a domain with finite measure then there exists a unique  minimizer $v_p$  in (\ref{rayplapapp}) satisfying the
normalizing conditions $v_p>0$ and $\int_{\Omega}v_p^pdx=1$. See the classical paper \cite{lqvsimple} and  also \cite{frala} for further discussions.  

By using the same argument of the proof of  Theorem~\ref{diff} combined with the results in \cite{lqv} concerning the continuous depencence of $v_p$ on $p$  (we  refer in particular to the local convergence result \cite[Thm.~6.3]{lqv} which by \cite{lieb} admits a natural global
version in the case of  domains of class $C^{1,\beta}$) one can prove the following theorem.

\begin{thm} Let $\Omega $ be a bounded domain in $\R^n$ of class $C^{1,\beta}$ with $\beta\in ]0,1]$. Then the function $p\mapsto \lambda_p$ is  
of class $C^1$ on $]1,\infty [$ and 
\begin{equation}\label{diffplap}
\lambda_p'=\int_{\Omega}|\nabla v_p|^p\ln |\nabla v_p|dx-\lambda_p\int_{\Omega}v_p^p\ln v_pdx \, , \quad
p\in]1,\infty [\, .
\end{equation}
\end{thm}

{\bf Acknowledgments.}  We  acknowledge financial support from 
the research project `Singular perturbation problems for differential operators' Progetto di Ateneo of the University of Padova.
The second author acknowledges financial support also from the research project  `INdAM GNAMPA Project 2015 - Un
approccio funzionale analitico per problemi di perturbazione singolare e di
omogeneizzazione'.
The second author is also member of the Gruppo Nazionale per l'Analisi Matematica, la Probabilit\`{a} e le loro Applicazioni (GNAMPA) of the
Istituto Nazionale di Alta Matematica (INdAM).

Both authors acknowledge the warm hospitality received by each other's institution on the occasion of several research visits.

\end{document}